\documentclass[11pt,a4paper]{article}

\usepackage[latin1]{inputenc}
\usepackage{verbatim}
\usepackage{enumerate}
\usepackage{array}
\usepackage{amsmath}
\usepackage{amssymb}
\usepackage{amsfonts}
\usepackage{theorem}
\usepackage{latexsym}
\usepackage{color}

\numberwithin{equation}{section}
\newcommand\NN{\mathbb{N}}

\newcommand\CC{\mathbb{C}}

\newcommand{\ve}{\varepsilon}
\newcommand{\no}{\noindent}
\newcommand{\beq}{\begin{equation}}
\newcommand{\eeq}{\end{equation}}

\begin{document}
\title{\bf{Distribution of points of interpolation and of zeros  of exact maximally convergent
multipoint Pad{\'e} approximants }}

\author{ R. K. Kovacheva,\\Institute of Mathematics and Informatics,
\\Bulgarian Academy of Sciences\\Sofia, Bulgaria}


\maketitle

{\bf Abstract:} Given a regular compact set $E$ in  $\CC$, a unit
measure $\mu$ supported by $\partial E,$  a triangular point set
 $\beta := \{\{\beta_{n,k}\}_{k=1}^n\}_{n=1}^{\infty},\beta\subset \partial E$
  and a function $f$, holomorphic  on $E$,
 let $\pi_{n,m}^{\beta,f}$ be the associated multipoint $\beta-$ Pad{\'e}
 approximant of order $(n,m)$.
  We show that if the sequence $\pi_{n,m}^{\beta,f}, n\in\Lambda, m-$ fixed,
  converges exact maximally to
 $f$, as $n\to\infty,n\in\Lambda$ inside the maximal  domain of
  $m-$ meromorphic continuability of $f$ relatively to the measure $\mu,$
  then the points $\beta_{n,k}$ are uniformly distributed  on $\partial E$
  with respect to the measure $\mu$ as $ n\in\Lambda$.
  Furthermore, a result about the zeros behavior of  the exact maximally
  convergent sequence $\Lambda$ is provided, under the condition that  $\Lambda$ is "dense enough."

{\bf Keywords:  } Multipoint Pad{\'e} approximants, maximal
convergence

{\bf MSC2010: }41A20, 41A21, 30E10
\section{Introduction}

We first introduce   some needed notations.

Let $\Pi_n,\, n\in\NN$ be  the class of the polynomials of degree
$\leq n$ and ${\cal R}_{n,m}:=\{r = p/q, p\in\Pi_n, q\in \Pi_m,\,
q\not\equiv 0\}$.

Given a compact set $E,$  we say that $E$ is {\it regular}, if the unbounded component of the complement $E^c:=\overline\CC\setminus E$  is solvable with respect to Dirichlet problem. We will assume throughout the paper that $E$ possesses a connected complement  $E^c$.
In what follows, we will be working with the $\hbox{max}-$ norm $||...||_E$ on $E$; that is $||...||_E:= \max_{z\in E}|...|(z).$

Let ${\cal B}(E)$ be  the class of the unit measures supported on  $E$; that is
$\hbox{supp}(...)\subseteq E.$ 
We say that the infinite sequence of Borel measures $\{\mu_n\}\in
{\cal B}(E)$ converges  in the weak topology to a measure $\mu$
and write $\mu_n\longrightarrow  \mu$,
 { if
 $$\int g(t) d\mu_n\to\int g(t) d\mu$$ for every function $g$ continuous  on $E$}.
 We associate with a measure $\mu\in{\cal B}(E)$ the logarithmic potential $U^\mu (z)$; that is,
 $$U^\mu (z):= \int\log\frac{1}{|z-t|}d\mu.$$
 Recall that   $U^\mu$ (\cite{SaffTotik}) is a function superharmonic in $\CC,$ subharmonic in $\overline \CC\setminus \hbox{supp}(\mu),$
  harmonic in $\CC\setminus \hbox{supp}(\mu) $ and
 $$U^\mu(z) = \ln \frac{1}{|z|} + o(1), z\to\infty.$$
 We also note the following basic fact (\cite{Carleson}):

  \no{\bf Carleson's lemma}: {\it Given the measures $\mu_1, \mu_2$ supported by $\partial E$, suppose that $U^{\mu_1}(z) = U^{\mu_2}(z)$ for every $z\not\in E$. Then $\mu_1 = \mu_2.$}

Finally, we associate with  a polynomial $p\in \Pi_n$  the normalized counting measure $\mu_p$ of $p, $ that is
$$\mu_p (F):= \frac{\hbox{number}\,\,\hbox{of}\,\,\hbox{zeros}\,\,\hbox{of}\, \,p\, \hbox{on}\,\, F}{\hbox{deg}\, p},$$
where $F$ is a point set in $\CC.$

Given a domain $B\subset \CC,$  a function $g$ and a number $m\in\NN$, we say that $g$ is $m-$meromorphic in $B$ ($g\in{\cal M}_m(B))$ if it has
no more than $m$  poles in $B$ (poles are counted with their multiplicities).
 We say that a function $f$ is holomorphic on the compactum $E$ and write $f\in{\cal A}(E)$, if it is holomorphic in some open neighborhood of $E$.


Let $\beta$ be  an infinite triangular  table of points, $\beta:=\{\{\beta_{n,k}\}_{k=1}^{n}\}_{ n=1,2,...}, \beta_{n,k}\in E$, with no limit points outside $E.$ Set $$\omega_n(z):=\prod_{k=1}^n(z-\beta_{n,k}).$$

 Let $f\in {\cal A}(E)$ and $(n,m)$ be a fixed pair of nonnegative integers.
The rational function $ \pi_{n,m}^{\beta,f} := p/q, $
 where the polynomials
$p\in\Pi_n$ and $q\in\Pi_m$ are such that
$$\frac{fq-p}{\omega_{n+m+1}}\in {\cal A}(E)$$ is called {\it a $\beta$-multipoint Pad{\'e} approximant of $f$ of order $(n,m)$ }.
As is well known, the function $\pi_{n,m}^{\beta,f}$ always exists and is unique (\cite{saff},\cite{Go2}). In the particular case when $\beta\equiv {0},$ the multipoint Pad{\'e} approximant $\pi_{n,m}^{\beta,f}$ coincides with the classical Pad{\'e} approximant  $\pi_{n,m}^{f}$ of order $(n,m)$ (\cite{Pe}).

{ Set
$$\pi_{n,m}^{\beta,f}:=\frac{P_{n,m}^{\beta,f}}{Q_{n,m}^{\beta,f}},\eqno({1})$$
where the polynomials $P_{n,m}^{\beta,f}$ and $Q_{n,m}^{\beta,f}$ do not have common divisors.
The zeros of $Q_{n,m}^{\beta,f}$ are called {\it free zeros} of $\pi_{n,m}^{\beta,f}$;   $\hbox{deg}\, Q_{n,m}\leq m$.}

We say that the points $\beta_{n,k}$ are {\it uniformly distributed relatively  to the  measure $\mu$,} if
$$\mu_{\omega_n}\longrightarrow \mu,\, n\to\infty.$$

We recall the notion of $m_1-$ Hausdorff measure
 (cf. \cite{Go3}).
 For $\Omega\subset\CC$, we set
    \begin{gather*}
    m_1(\Omega) := {\rm inf} \left\{\sum_\nu |V_\nu| \right\}
    \end{gather*}
where the infimum is taken over all coverings $\{V_\nu\}$ of $\Omega$
by disks  and $|V_\nu|$ is the radius of the disk $V_\nu$.

Let $D$ be a domain in $\CC$ and $\varphi$ a function defined in
$D$ with values in $\overline{\CC}$. A sequence of functions
$\{\varphi_n\}$, meromorphic in $D$, is said to converge to a
function $\varphi$ {\it $m_1$-almost uniformly inside $D$} if for any compact subset $K
\subset D$ and every  $\varepsilon > 0$ there exists a set
$K_\varepsilon \subset K$ such that $m_1(K \setminus
K_\varepsilon) < \varepsilon$ and the sequence $\{\varphi_n\}$
converges uniformly to $\varphi$ on $K_\varepsilon$.

 For $\mu\in{\cal B}(E),$ define
 $$\rho_{\hbox{min}}:= \inf_{z\in E} e^{- U^\mu(z)}$$
  and $$\varrho_{\hbox{max}}:= \max_{z\in E} e^{- U^\mu(z)};$$
  ($U^\mu$ is superharmonic on $E$; hence it attains its minimum (on $E$)).
   As is known (\cite{Ts}, \cite{SaffTotik}), $$e^{-U^\mu(z)}\geq \rho_{\hbox{min}},\, z\in E^c.$$
    Set,  for $r >  \rho_{\hbox{min}},$
 $$E_\mu(r):= \{z\in\CC, e^{-U^\mu(z)} < r\}.$$
 Because of the upper  semicontinuity of the function $\chi(z):=e^{- U^\mu(z)}$, the set $E_\mu(r)$  is open; clearly $E_\mu(r_1)\subset E_\mu(r_2)$ if $r_1\leq r_2$
and  $E_\mu(r) \supset E$ if $r >\varrho_{\hbox{max}}.$

 Let $f\in {\cal A}(E)$ and $m\in\NN$ be fixed. Let $R_{m,\mu}(f)=R_{m,\mu}$ and $D_{m,\mu}(f)=D_{m,\mu}:=E_\mu(R_{m,\mu})$ denote, respectively,
 {\it the radius and domain of $m-$ meromorphy with respect to $\mu$;} that is
 $$R_{m,\mu}:= \sup\{r, f\in {\cal M}_m(E_\mu(r))\}.$$

 Furthermore, we  introduce the notion  of a {\it  $\mu-$maximal convergence  to  $f$ with respect to the $m-$meromorphy of a sequence of rational functions
 $\{r_{n,\nu}\}$
} (a $\mu$-maximal convergence): that is,
 for any $\ve >0$ and each compact set $K\subset D_m$, there exists a set $K_\ve\subset K$ such that $m_1(K\setminus K_\ve)<\ve$ and
$$\limsup_{n+\nu\to\infty}||f-r_{n,\nu}||_{K_\ve}^{1/n}\leq \frac{||e^{-U^\mu}||_K}{R_{m,\mu}(f)}.$$

  Hernandez and Calle Ysern proved the following:

  \no{\bf Theorem A, \cite{spain}} :
  {\it Let $E, \mu, \beta$ and $\omega_n,n=1,2,...,$ be defined as above. Suppose that $\mu_{\omega_n}\longrightarrow \mu$ as $n\to\infty$ and $f\in {\cal A}(E).$ Then,  for each fixed $m\in\NN$,  the sequence $\pi_{n,m}^{\beta, f}$ converges to $f$
  $\mu-$maximally with respect to the $m-$ meromorphy.
}

Theorem A generalizes E. B. Saff's  theorem of Montessus de Ballore's type about multipoint Pad{\'e} approximants (see\cite{saff}).

 We now utilize the {\it normalization of the polynomials $Q_{n,m}(z)$
 with respect to a given open set $D_{m,\mu};$}  that is,
$$Q_{n,m}(z) = \prod(z-\alpha'_{n,k})\prod (1-z/\alpha''_{n,k}),\eqno{(2)}$$ where $\alpha'_{n,k},\,\alpha''_{n,k}$ are the zeros lying inside, resp. outside $D_{m,\mu}$.
Under this normalization, for every compact set $K$ and $n$ large
enough there holds
$$\Vert Q_{n,m}^{\beta, f}\Vert_K\leq C_1,$$
where $C_1 = C_1(K)$ is a positive constant, depending on $K.$
In the sequel, we denote by $C_i$ positive constant, independent on $n$ and different at different occurances.

In \cite{spain}, the set $K_\ve$ (look at the definition of a $\mu-$maximal  convergence) is explicitly written, namely$K_\ve:= K\setminus \Omega(\ve),$ where
 $$\Omega(\ve):=\bigcup_{n=1}^\infty(\bigcup_{\alpha_{n,k}'}\{z, |z-\alpha_{n,k}'| < \ve/(2mn^2)\}).$$
For $\Omega(\ve)$ we have
$$m_1(\Omega(\ve))\leq \ve.$$
For points $z\not\in \Omega(\ve)$, we have
$$|Q_{n,m}^{\beta, f}(z)| \geq C_2(\ve/mn^2)^{k_n},$$
where $k_n$ stands for the number of the zeros of $Q_{n,m}^{\beta, f}$ in $D_{m,\mu}; k_n\leq m.$

Let $Q$ be the monic polynomial, the zeros of which coincide with
the poles of $f$ in $D_{m,\mu}; \,\hbox{deg}\, Q\leq m.$ It was
proved in \cite{spain} (Proof of Lemma 2.3)
 that for every compact
subset $K$ of $D_{m,\mu}$
$$\limsup_{n\to\infty}\|fQQ_{n,m}^{\beta, f}-QP_{n,m}^{\beta, f}\|^{1/n}_K\leq  \frac{||e^{-U^\mu}||_K}{R_{m,\mu}}.\eqno{(3)}$$
 Hence, the function $-U^\mu(z) - \ln R_{m,\mu}$ is a harmonic majorant in $D_{m,\mu}$
 of the family $\{|(fQQ_{n,m}^{\beta, f}-QP_{n,m}^{\beta, f})(z)|^{1/n}\}_{n=1}^\infty.$

\no{\bf Theorem B, \cite{spain}} {\it With  $E, \mu,\, m,\omega_n$ and $ f$
as in Theorem A, assume that $K$ is a regular compact set  for which $\|e^{-U^\mu}\|_K$ is not attained at a  point on $E$.
Suppose that the function $f$ is defined on $K$ and satisfies
 $$\limsup_{n\to\infty} \Vert f - \pi_{n,m}^{\beta, f}\Vert_{K}^{1/n} \leq \Vert e^{-U^\mu}\Vert_K/R < 1.$$
Then }$R \leq R_{m,\mu}(f)
.$

\no{\bf Remark:} Suppose that $\infty > R_{m,\mu} > \varrho_{\hbox{max}}$ and $D_{m,\mu}$ is connected.
Let $V$ be a disk in $D_{m}\setminus E_\mu(\varrho_{\hbox{max}}),$ centered at a point $z_0$ of radius $r>0$ and such that $f$ is analytic on $V.$ Fix  $r_1,\,0 < r_1 <r $ and set $A:=\{z, r_1 \leq |z-z_0| \leq r\}.$
Fix  a number $\ve < (r-r_1)/4.$ Introduce, as before, the set  $\Omega(\ve)$.
Recall that
$$m_1(\Omega(\ve))\leq \ve.$$
It is clear that the set $A\setminus \Omega(\ve)$ contains a
concentric circle $\Gamma$ (otherwise we would obtain  a
contradiction with $m_1(\Omega(\ve)) < (r-r_1)/4.$) We note that
the function $f$ and the rational functions $\pi_{n,m}^{\beta,f}$
are well defined on $\Gamma.$ Viewing  (3), we may write
$$\limsup_{n\to\infty} \|QQ_{n,m}^{\beta, f}f - QP_{n,m}^{\beta, f}\|^{1/n}_\Gamma \leq \|e^{-U^\mu}\|_\Gamma/R_{m,\mu}.$$
 Suppose that
 $$\limsup_{n\to\infty} \|QQ_{n,m}^{\beta, f}f-QP_{n,m}^{\beta, f}\|^{1/n}_\Gamma < \|e^{-U^\mu}\|_\Gamma/R_{m,\mu},$$
  or, what is the same,
  $$\limsup_{n\to\infty} \|QQ_{n,m}f-QP_{n,m}^{\beta, f}\|^{1/n}_\Gamma  \leq
    \|e^{-U^\mu}\|_\Gamma/(R_{m,\mu}+\sigma) < 1$$ for an appropriate $\sigma > 0. $
Then, $$ \vert (f - \pi_{n,m}^{\beta,f})(z)\vert_{\Gamma} \leq C_3  (n^2m/\ve)^m(\|e^{-U^\mu}\|_\Gamma/(R_{m,\mu}+\sigma))^n $$
for all $z\in\Gamma$ and  $n$ large enough. This leads to
$$\limsup_{n\to\infty}\Vert f - \pi_{n,m}^{\beta,f}\Vert_{\Gamma}^{1/n}\leq \|e^{-U^\mu}\|_\Gamma/(R_{m,\mu}+\sigma).$$
using  Theorem B, we arrive at $R_{m,\mu}+\sigma < R_{m,\mu}.$ The
contradiction yields
 $$\limsup_{n\to\infty} \|QQ_{n,m}^{\beta, f}f-QP_{n,m}^{\beta, f}\|^{1/n}_{\overline V_\Gamma} = \|e^{-U^\mu}\|_{\overline V_\Gamma}/R_{m,\mu},
$$
where $V_\Gamma$ is the disk bounded by  $\Gamma.$

Then the function $-U^\mu-\ln R_{m,\mu}$ is an exact harmonic
majorant of the family $\{\vert fQQ_{n,m}^{\beta, f}-QP_{n,m}^{\beta, f}\vert^{1/n}\}$
in $D_{m,\mu}$ (see  (3)). Therefore,   there exists a subsequence
$\Lambda$ such that for every compact subset $K\subset
D_{m,\mu}\setminus E $
$$\lim_{n\in\Lambda} \|QfQ_{n,m}^{\beta, f}-P_{n,m}^{\beta, f}Q\|_K^{1/n} =  \Vert e^{-U^\mu}\Vert_K/R_{m,\mu}.\eqno{(4)}$$
(see {\cite{Wa1},\cite{Wa2}) for a discussion of exact harmonic majorant)).
We will refer to this sequences as to {\it  an exact maximally convergent sequence.}

We prove

\no{\bf Theorem 1:} {\it Under the same conditions on $E$, assume
that $\mu\in {\cal B}(\partial E)$ and   that $\beta\subset
\partial E$ is a triangular  set of points. Let $m\in \NN$ be
fixed, $f\in {\cal A}(E)$ and $\varrho_{\hbox{max}} < R_{m,\mu} <
\infty$. Suppose that $D_{m,\mu}$ is connected. If for a
subsequence $\Lambda$ of the multipoint Pad{\'e} approximants
$\pi_{n,m}^{\beta, f}$ condition (4) holds, then
$\mu_{\omega_n}\longrightarrow \mu$ as
 $n\to\infty,\, n\in\Lambda$.}

The problem of the distribution of the points of interpolation of
multipoint Pad{\'e} approximants was investigated, so far,  only
for the case when the measure $\mu$ coincides with the equilibrium
measure $\mu_E$ of the compact set $E.$ It was first rased by J.
L. Walsh (\cite{Wa3}, Chp. 3) while considering maximally
convergent polynomials with respect to the equilibrium measure. He
showed that the sequence $\mu_{\omega_n}$ converges weakly to
$\mu_E$ through the entire set $\NN$ (respectivey their associated
measures onto the boundary of $E$) iff the interpolating
polynomials of every function $f_t(z)$ of the form $f_t(z):=
(t-z), t \not\in E, t -$ fixed, converge maximally to $f_t.$
Walsh's result was extended to multipoint Pad{\'e} approximants
with a fixed number of the free poles by N. Ikonomov in
\cite{Ikonomov}, as well as to generalized Pad{\'e} generalized
approximants, associated with a regular condenser
(\cite{ikonomov}). { The
 case of polynomial interpolation of an arbitrary function $f$ holomorphic in $E$  was
considered by R. Grothmann (\cite{grothman}); he established  the
existence of an appropriate sequence $\Lambda$ such that
$\mu_{\omega_n}\longrightarrow \mu_E,\, n\to\infty,\,n\in\Lambda
,$ respectively the balayage measures onto $\partial E.$}
Grothmann's result was generalized in relation to multipoint
Pad{\'e} approximants $\pi_{n,m}^{\beta, f}$ with a fixed number
of the free poles (see \cite{sps}). { Finally, in \cite{BlKov} was
considered the case when the degrees of the denominators tend
slowly to infinity, namely $m_n = o(n/\ln n).$}

 As a consequence of Theorem 1, we derive

\no{\bf Theorem 2:} {\it Under the conditions of Theorem 1, suppose that the exact maximally convergent  sequence $\Lambda :=\{n_k\}_{k=1}^
\infty$
satisfies the condition to be "dense enough"; that is
$$\limsup\frac{n_{k+1}}{n_k} < \infty.$$ Then there
is at least one point $z_0\in\partial D_{m,\mu}(f)$ such that for
every  disk $V_{z_0}(r)$ centered at $z_0$ of  radius $r$
 }  $$ \limsup_{n\to\infty,\,
n\in\Lambda}\mu_{P_{n,m}^{\beta, f}}(V_{z_0} (r)) >
0,\,\hbox{as}\, n\to\infty,  n\in\Lambda.$$

 \no{\bf Proof of Theorem 1:} Set $Q_{n,m}^{\beta,f}:=Q_n,\, P_{n,m}^{\beta,f}:= P_n$ and $F:=fQ.$
 Fix numbers $R,\tau, r$ such that $ \varrho_{\hbox{max}} < R <\tau <  r < R_{m,\mu}$ and $E_\mu(R)$ is connected. Then, by the conditions of the theorem,
  for every compactum $K\subset D_{m,\mu}$ (comp.(4))
 $$\lim_{n\in\Lambda} \|FQ_{n}-QP_{n}\|_K^{1/n} =  \Vert e^{-U^\mu}\Vert_K/R_{m,\mu},n\in\Lambda.\eqno{(5)}$$
Select a positive number $\eta$ such that $R+\eta < \tau < \tau+\eta<r
< R_{m,\mu}.$ Let $\Gamma$  be an analytic curve in
$E_\mu(r)\setminus E_\mu({\tau+\eta})$ such that $\Gamma$ winds around every point in $E_\mu(\tau)$ exactly once.
In an analogous way, we select a curve $\gamma\subset
E_\mu({R+\eta})\setminus E_\mu(R).$ Additionally, we require that $U^\mu$ is constant on $\Gamma$ and $\gamma.$
 Set $$F_n(z):= \frac{1}{n}\ln
|FQ_{n}-P_nQ|(z) + U^\mu(z) + \ln R_{m,\mu},n\in\Lambda.\eqno{(6)}$$ Let $\sigma > 0$ be arbitrary. The functions $F_n$
are subharmonic in $E_\mu(r)\setminus E_\mu(R).$  By (5) and the choice of $\Gamma,$
$$\max_{t\in \Gamma}F_n(t)\leq -\min_{t\in \Gamma} +\max_{t\in \Gamma} +\sigma \leq \sigma, N\in\Lambda,
 n \geq n_1 0 n_1(\sigma)$$
and, analogously,
$$\max_{t\in \gamma}F_n(t)\leq -\min_{t\in \Gamma} +\max_{t\in \Gamma} \leq \sigma, N\in\Lambda, n>n_1
.$$ Then, by the max-principle of subharmonic functions,
$$\max_{z\in {\cal A}_{\gamma, \Gamma}}F_n(z)\leq \sigma, n\in\Lambda, n\geq n_1, N\in\Lambda,\eqno{(7)} $$
where ${\cal A}_{\gamma, \Gamma}$ is the "annulus",  bounded by $\Gamma$ and $\gamma.$

On the other hand, by  (5), for any compact set $K\subset
E_{r}\setminus E_R$ and $n$ large enough there is a point
$z_{n,K}\in K$ such that
$$-\min_K U^\mu (z_{n,K}) - \ln R_{m,\mu}-\sigma
\leq \frac{1}{n}\ln|FQ_n(z_{n,K})-QP_n(z_{n,K})|, n\geq
n_3(K),n\in\Lambda.$$ Therefore, $$-\sigma \leq F_n(z_{n,K}),
n\geq n_2(K,\sigma). \eqno{(8)}$$ Further, by the formula of
Hermite-Lagrange, for $z\in\gamma$ we have
$$FQ_n(z)-QP_n(z) =
\frac{1}{2\pi
i}\int_{\Gamma}\frac{\omega_{n+m+1}(z)}{\omega_{n+m+1}(t)}\frac{FQ_n(t)-QP_n(t)}{t-z}dt.$$
Hence, by (5), $$\frac{1}{n}\ln\vert FQ_n(z)-QP_n(z)\vert \leq $$
$$\max_{t\in\Gamma}
U^{\omega_{n+m+1}}(t) - U^{\omega_{n+m+1}}(z) +
\frac{1}{n}\ln ||FQ_n-QP_n||_\Gamma + \frac{1}{n}
\hbox{const} \leq $$
$$\max_{t\in\Gamma}
U^{\omega_{n+m+1}}(t) - U^{\omega_{n+m+1}}(z) - \min_{t\in\Gamma}
U^\mu (t) - \ln R_{m,\mu} + \sigma,n\in\Lambda, n\geq n_3 =
n_3(\sigma) > n_1,$$ where $U^{\omega_{n+m+1}}:=
U^{\mu_{\omega_{n+m+1}}}.$ To simplify the notations, we set
$U^{\omega_{n+m+1}}:= U^{\omega_{n}}.$ ( The correctness will be
not lost,  since $m\in\NN$ is fixed). Involving into consideration
the functions $F_n$ (see (6)),
 we get for
$z\in\gamma$
$$ F_n(z) \leq \max_{t\in\Gamma}(U^{\omega_n}(t)
- U^\mu (t)) + \max_{t\in\Gamma} U^\mu (t) +$$
$$(U^\mu (z) - U^{\omega_{n}}(z)) - \min U^\mu (t) +\sigma, n\in\Lambda, n \geq
n_2\geq n_1.$$ By Helly's selection theorem (\cite{SaffTotik} ),
there exists a subsequence of $\Lambda$ which we denote again by
$\Lambda$ such that
$\mu_{\omega_{n+m+1}}:=\mu_{\omega_n}\longrightarrow \omega,\,
n\in\Lambda.$ Passing to the limit,  we obtain
$$\limsup_{\Lambda}\vert F_n(z)\vert \leq \max_{t\in\Gamma}(U^{\omega}(t) - U^\mu (t))
+(U^\mu (z) - U^{\omega}(z)), z\in\gamma.
\eqno{(9)}$$

Consider  the function $\phi$, harmonic
in ${\cal A}_{\Gamma, \gamma}$ and
$$\phi:=\left\{\begin{array}{ll}
0, &\Gamma,\\
 \min (0, -\min_{t\in\gamma}(U^\mu(t) - U^{\omega}
(t)) + ( U^\mu (z) - U^{\omega}(z)
),&\gamma\\
\end{array}
\right.
$$
From (7) and (9), we arrive at
$$\limsup  F_n(z) \leq \phi,$$
for $z$ in  ${\cal A}_{\Gamma,\gamma}$.
Being harmonic, $\phi$ obeys the maximum and the minimum principles in this region. The definition yields
\[ \phi (z) \leq 0, z\in {\cal {A}}_{\Gamma,\gamma} \]
We will show that
$$
\phi (z) \equiv 0, \eqno{(10)}
$$

 Suppose that (10) is not true. Let $\Upsilon$ be a closed curve in the
  set $E_{R+\eta}-\gamma^o,$
 where $\gamma^o$ stands for the interior of $\gamma.$
Then there exists a number  $\theta > 0$ such  that $\phi \leq -\Theta$ for every $z\in \Upsilon$.
 This inequality contradicts (8), for $\sigma$ close enough to the zero and $n\in\Lambda$ sufficiently large. .

Hence, $\phi\equiv 0$. Then  the definition of $\phi$ yields
\[ U^{\mu}(z) - U^{\omega}(z) \equiv  \min_{t\in \gamma} \big(U^{\mu}(t) -  U^{\omega}(t)\big) ,\ z\in\gamma. \]
The function $U^{\mu}(z) - U^{\omega}(z)$ is harmonic in the unbounded complement $G$ of $\gamma$,  and by the maximum principle,
\[ U^{\mu}(z) - U^{\omega}(z) \equiv  \hbox{Constant,} \ z\in G; \]
 consequently,
\[ U^{\mu}(z) - U^{\omega}(z) \equiv  \hbox{Constant,} \ z\in E^c. \]
On the other hand, $(U^{\mu} - U^{\omega})(\infty)= 0$, which yields
$U^\mu\equiv U^\omega$ in $E^c.$
By Carleson's Lemma, $\mu=\omega.$ On this, Theorem 1 is proved.\hfill{\bf Q.E.D.}

\medskip
The proof of Theorem 2 will be preceded by an auxiliary lemma

\no{\bf Lemma 1, \cite{rkk}  :} {\it Given a domain $U$,
 a regular compact subset $S$ and a sequence $\vartheta := \{n_k\}$
  of positive integers, $n_k < n_{k+1},\, k = 1,2,\cdots  ,$
  such that $$\limsup \frac{n_{k+1}}{n_k} < \infty.$$
  Suppose that $\{\phi_{n_k}\}$ is a sequence of rational functions,
   $\phi_{n_k}\in{\cal R}_{n_k,n_k},\, k - 1,2,\cdots ,\phi_{n_k} = \phi'_{n_k}/\phi"_{n_k}$
   having no more that $m$ poles in $U$ and converging
   uniformly of $\partial S$ to a function $\phi\not\equiv 0$ such that $$\limsup_{n_k\to\infty, n_k\in\Lambda}\|\phi_{n_k}- \phi\|_{\partial S}^{1/n} < 1.$$
Assume, in addition, that on each compact subset of $U$
$$\lim\mu_{\phi_{n_k}'}(K) \longrightarrow 0.\eqno{(11)}$$

Then the function $\phi$ admits a continuation into $U$ as a meromorphic function with no more that $m$ poles. }

\smallskip

\no{\bf Proof of Theorem 2}. We preserve the notations in the
proof of  Theorem 1.

The proof of Theorem 2   follows from Lemma 1 and Theorem 1.
 Indeed, under  the conditions of the theorem
 the sequence $\{\pi_{n}\}_{n\in\Lambda}$
 converges maximally to $f$ with
  respect to the measure $\mu$ and the domain $D_{m,\mu}.$
   Hence, inside  $D_{m,\mu}$ condition (11) if fulfilled.
   From the proof of Theorem 1, we see that  there is a regular compact subset $S$ of $D_{m,\mu}$
    such that $\limsup_{n\in\Delta}\|f-\pi_n\|_S^{1/n} < 1.$

   Suppose now that the statement of Theorem 2 is not true.
   Then there  is an open  { strip $W$} containing  $\partial D_{m,\mu}$
    such that on each compact subset
    of $W$ condition $(11)$ holds.
     Applying Lemma 1 with respect to the sequence $\pi_n$ and the domain $D_{m,\mu}\bigcup W,$
      we conclude   that $f\in{\cal M}_m(\overline{D_{m,\mu}}).$
      This contradicts  the definition of $D_{m,\mu}.$

On this, the proof of Theorem 2 is completed. \hfill{\bf Q.E.D.}

\smallskip
Using again Lemma 1 and applying Theorem A, we obtain a result related to the zero distribution of the  sequence $\{\pi_{n,m}^{\beta, f}\}.$

\smallskip
\no{\bf Theorem 3:} {\it Let $E$ be a regular compactum in $\CC$
with a connected complement, let $\mu\in{\cal B}(E)$ and $\beta\in
E$ be a triangular point set.
 Let the polynomials $\omega_n,n=1,2,...,$ be defined as above.
  Suppose that $\mu_{\omega_n}\longrightarrow \mu$ as $n\to\infty$ and $f\in {\cal A}(E).$
  Let $m\in\NN$ be fixed, and suppose that   
  $R_{m,\mu} < \infty.$ Then,  there is at least one point $z_0\in\partial D_{m,\mu}$ such that
  $\limsup_{n\to\infty} \mu_{\pi_{n,m}^{\beta, f}}
  (\overline V_{z_0}(r)) > 0$ for every positive $r$.}

Ralitza K. Kovacheva\\Institute of Mathematics and Informatics, \\Bulgarian Academy of Sciences\\Acad. Bonchev str. 8, \\1113 Sofia, Bulgaria\\rkovach@math.bas.bg

\end{document}